\documentclass[11pt,a4paper,fleqn]{article}
\pdfoutput=1
\RequirePackage{amsthm,amsmath,natbib}
\RequirePackage{amssymb,amstext}
\RequirePackage{hypernat}
\usepackage[pdfpagemode=UseOutlines ,plainpages=false
,hypertexnames=false ,pdfpagelabels ,hyperindex=true,colorlinks=true]{hyperref}
	\makeatletter%
	\Hy@breaklinkstrue%
	\makeatother%
\usepackage{color}
\definecolor{darkred}{rgb}{0.6,0.0,0.1}
\definecolor{darkgreen}{rgb}{0,0.5,0}
\definecolor{darkblue}{rgb}{0,0,0.5}
\hypersetup{colorlinks ,linkcolor=darkblue ,filecolor=darkgreen
,urlcolor=darkblue ,citecolor=black}

\usepackage{		url}
\usepackage{verbatim}
\usepackage{ bm}

\renewcommand{\cite}{\citet}
\usepackage{jan-abrev-package}

\definecolor{dgreen}{rgb}{0,0.5,0}
\definecolor{dblue}{rgb}{0,0,0.9}
\definecolor{dred}{rgb}{0.6,0.0,0.1}
\definecolor{dgold}{rgb}{0.5,0.3,0.0}
\definecolor{dvio}{rgb}{0.6,0.3,0.5}
\definecolor{gray}{rgb}{0.5,0.5,0.5}

\newcommand{\dr}{\color{dred}}

\newtheoremstyle{mysc}
  {3pt}
  {3pt}
  {\it}
  {}
  {\color{darkred}\sc}
  {.}
  {.5em}
  {}

\newtheoremstyle{myex}
  {10pt}
  {10pt}
  {\rm}
  {}
  {\color{darkred}\sc}
  {.}
  {.5em}
  {}

\theoremstyle{mysc}\newtheorem{prop}{Proposition}[section]
\theoremstyle{mysc}\newtheorem{assumption}{Assumption}[section]
\theoremstyle{mysc}
\theoremstyle{mysc}\newtheorem{theo}[prop]{Theorem}
\theoremstyle{mysc}
\theoremstyle{mysc}\newtheorem{lem}[prop]{Lemma}
\theoremstyle{myex}\newtheorem{rem}{Remark}[section]
\theoremstyle{myex}
\theoremstyle{myex}

\numberwithin{equation}{section}

\author{{\sc Jan Johannes}\thanks{Universit\"at Heidelberg, Institut f\"ur Angewandte Mathematik, Im Neuenheimer Feld, 294, D-69120 Heidelberg, Germany, e-mail: \url{johannes@statlab.uni-heidelberg.de}}}

\title{{\bf Nonparametric estimation in functional linear models with second order stationary regressors.}}

\begin{document}
\date{\today}
\maketitle


\begin{abstract} We consider the problem of estimating the slope parameter in  functional linear regression, where scalar responses
$Y_1,\dotsc,Y_n$ are modeled in dependence of second order stationary random functions
$X_1,\dotsc,X_n$.  An orthogonal series estimator of the functional slope parameter with additional thresholding in the Fourier domain  is proposed and its performance is measured  with respect to a wide range of weighted risks covering as examples the mean squared prediction error  and the mean integrated squared error for derivative estimation.  In this paper the minimax optimal rate of convergence of the estimator is derived over a large class of different regularity spaces for the slope parameter and of different link conditions for the covariance operator.  These general results are illustrated by the particular example of the well-known Sobolev space of periodic functions as regularity space for the slope parameter and the case of finitely or infinitely smoothing  covariance operator.
\end{abstract}

\begin{tabbing}
\> Mean squared error of prediction, Minimax theory, Sobolev space.\\[.2ex]
\end{tabbing}

\section{Introduction}\label{sec:intro}
Functional linear models have become very important in a diverse range of disciplines,
including medicine, linguistics, chemometrics as well as econometrics (see for instance 
\cite{RamsaySilverman2005} and \cite{FerratyVieu2006}, for several case studies, or more specific, \cite{ForniReichlin1998} and \cite{PredaSaporta2005} for  applications in economics).   Roughly speaking, in all
these applications the dependence 
of a response variable $Y$  on the variation of an explanatory random function $X$ is modeled by 
\begin{equation}\label{intro:e1}Y=\int_{0}^1\beta(t)X(t)dt +\sigma\epsilon,\quad\sigma>0,
\end{equation}
for some error term $\epsilon$. One objective is then to estimate nonparametrically the
slope function $\beta$  based on an independent and identically distributed (i.i.d.) sample of $(Y,X)$.  

In this paper we suppose  that the random function $X$ is taking  its  values in $L^2[0,1]$, which
is endowed with the usual inner product $\skalar$ and induced norm $\norm$, and that $X$ has a finite second moment, i.e., $\Ex \normV{X}^2<\infty$.  In order to simplify notations we assume that the mean function of $X$ is zero. Moreover, the random function $X$ and the error term $\epsilon$ are uncorrelated, where  $\epsilon$ is assumed to have  mean zero and variance one. This situation
has been considered, for example, in \cite{CardotFerratySarda2003} or \cite{MullerStadtmuller2005}. Then multiplying both sides in (\ref{intro:e1}) by $X(s)$ and taking the expectation leads to 
\begin{equation}\label{method:e1}g(s):=\Ex[YX(s)]=
\int_{0}^1\beta(t)\cov(X(t),X(s))dt=:[T_{\cov} \beta](s),\quad s\in[0,1],
\end{equation}
where $g$ belongs to $L^2[0,1]$ and $T_{\cov}$ denotes the covariance operator associated to the random function $X$.  
Estimation of $\beta$ is thus linked with the
inversion of the covariance operator $T_{\cov}$ of $X$ and, hence called an inverse problem. We assume that there exists a unique solution $\beta\in L^2[0,1]$ of equation (\ref{method:e1}),  i.e.,  $g$ belongs to the range $\cR(T_{\cov})$ of $T_{\cov}$, and $T_{\cov}$ is injective. However, as usual in the context of inverse problems all the results below could also be obtained straightforward for the unique  least-square solution with minimal norm, which exists if and only if $g$ is contained in the direct sum of $\cR(T_{\cov})$ and its orthogonal complement $\cR(T_{\cov})^\perp$ (for a definition and detailed discussion in the context of inverse problems see chapter 2.1 in \cite{EHN00}, while in the special case of a functional linear model we refer to \cite{CardotFerratySarda2003}). 

The normal equation (\ref{method:e1}) is the continuous equivalent of a normal equation \lq\lq$\Ex XY=\Ex XX^t \beta$\rq\rq\ in a linear model \lq\lq$Y=X^t\beta+\epsilon$\rq\rq, where the  covariance matrix \lq\lq$\Ex XX^t$\rq\rq\ has always a continuous generalized inverse. However, due to the finite second moment of $X$ the covariance operator $T_{\cov}$ of $X$ defined in (\ref{method:e1}) is nuclear (c.f. \cite{DauxoisPousseRomain82}). Thereby, unlike in the   linear model,  a continuous generalized inverse of $T_{\cov}$ does not exist if the range of the operator $T_{\cov}$ is an infinite dimensional subspace of $L^2[0,1]$. 
This corresponds to the setup of ill-posed inverse problems (with the additional
difficulty that $T_{\cov}$ in \eqref{method:e1} is unknown and hence, has to be estimated). 

In the literature several approaches are proposed in order to circumvent the instability issue due to an inversion of $T_{\cov}$.  
Essentially, all of them replace  the operator $T_{\cov}$ in equation (\ref{method:e1}) by a regularized version having a continuous generalized inverse. 
A popular  example is based on a functional principal
components  regression (c.f. \cite{Bosq2000}, \cite{CardotMasSarda2007} or \cite{MullerStadtmuller2005}), which corresponds to a  method  called spectral cut-off in the numerical analysis literature (c.f. \cite{Tautenhahn96}). An other example is the Tikhonov regularization (c.f. \cite{HallHorowitz2007}), where the regularized solution $\beta_\alpha$ is defined as unique minimizer of the Tikhonov functional 
$F_\alpha(\beta)=\normV{T_{\cov} \beta -g}^2 + \alpha \normV{\beta}^2$ for some strictly positive $\alpha$.  A regularization through
a penalized least squares approach after projection onto some basis (such as
splines) is also considered in \cite{RamsayDalzell1991}, \cite{EilersMarx1996} or \cite{CardotFerratySarda2003}.

In opposite to the model assumptions considered until now in the literature in this paper we suppose that the regressor $X$ is second order stationary. Over relatively short periods of time, the assumption of  second order stationarity is  in many situations realistic and can be checked from the data by estimating the covariance function using the multiple realizations of $X$.  Moreover, assuming second order stationarity allows us to generalize the known results in essentially two directions. First, we can unify the measures of performances for the estimator as considered in the literature  and second it is possible to present a simple estimation strategy which is optimal in a minimax sense over a wide range of possible regularity spaces for the slope functions $\beta$ as well as various forms of link conditions for the covariance operators $\op$. To be  more detailed:

In this paper we show that in case of second order stationary regressor $X$ the associated covariance operator $\op$ admits a spectral decomposition $\{\lambda_j,\psi_j,j\geq1\}$ given by the trigonometric basis $\{\psi_j\}$ (defined below) as eigenfunctions and a strictly positive, possibly not ordered, zero-sequence $\lambda:=(\lambda_j)_{j\geq 1}$ of corresponding eigenvalues.  Then the normal equation can be rewritten as follows 
\begin{equation}\label{intro:e2}
\beta=\sum_{j=1}^\infty \frac{g_j}{\lambda_j}\cdot\psi_j\quad\text{ with } {g}_j:=\skalarV{{g},\psi_j},\; j\geq1.
\end{equation} 
It is well-known that even in case of an a-priori known sequence $\lambda$ of eigenvalues    replacing in (\ref{intro:e2}) the unknown function $g$ by a consistent estimator $\widehat{g}$ does in general not lead  to a $L^2$-consistent estimator  of $\beta$. To be more precise, since  $\lambda$  is a zero-sequence, $\Ex\normV{\widehat{g}- g}^2=o(1)$ does generally not imply $\sum_{j=1}^\infty \lambda_j^{-2}\cdot\Ex|\skalarV{\widehat{g}-g,\psi_j}|^2 =o(1)$, i.e., the inverse operation of the covariance operator $T_{\cov}$ is not continuous. Essentially, all of the approaches mentioned above circumvent this instability issue  by replacing equation (\ref{intro:e2}) by a regularized version which avoids that the denominator becomes too small. For instance, in case of  a Tikhonov regularization (c.f. \cite{HallHorowitz2007}) in (\ref{intro:e2})  the factor $1/\lambda_j$ is replaced by $\lambda_j/(\alpha+\lambda_j^2)$.

In the literature so far the performance of an estimator of $\beta$ has been measured either by considering a squared prediction error or an integrated squared error. We show in this paper that these approaches can be unified  by considering a loss given by a weighted norm. To be more precise for $f \in L^2[0,1],$ we define
\begin{equation}\label{intro:e3}
\normV{f}_\hw^2 := \sum_{j=1}^\infty \hw_j |\skalarV{f,\psi_j}|^2
\end{equation}
for  some strictly positive sequence of weights  $\hw:=(\hw_j)_{j\geq 1}$. Then, the performance of an estimator $\widehat{\beta}$ of $\beta$  is measured by the $\cF_{\hw}$-risk, that is $\Ex\|\widehat{\beta}-\beta\|_{\hw}^2$.  This general framework allows us with an appropriate choice of the weight sequence $\hw$ to cover both, the risk in terms of mean integrated squared error, i.e., $\hw\equiv 1$,  as well as the mean squared prediction error. Indeed,   the squared prediction  error of a new value of $Y$ 
given any random function $X_{n+1}$ possessing the same distribution as $X$ and being
independent of $X_1,\dotsc,X_n$ can be evaluated as follows (see for example \cite{CardotFerratySarda2003} or \cite{CrambesKneipSarda2007} for similar setups)
\begin{equation*}
\Ex \Bigl[ \left| \skalarV{\widehat\beta,X_{n+1}} - \skalarV{\beta,X_{n+1}} \right|^2\, \Bigl\vert\, \widehat\beta\Bigr]
=  \skalarV{\op (\widehat\beta-\beta), (\widehat\beta-\beta)}=\sum_{j\geq 1}\lambda_j |\skalarV{\widehat\beta-\beta,\psi_j}|^2,
\end{equation*}
where we have used for the last identity that the regressor is second order stationary, i.e, $\op$ admits $\{\lambda_j,\psi_j,j\geq 1\}$ as spectral decomposition. Consequently, choosing $\hw \equiv  \lambda$ the $\cF_{\hw}$-risk is equivalent to the mean squared prediction error. We present this specific  situation in Section \ref{sec:ex}  below. It is worth to note, that 
the $L^2$-norm $\normV{f^{(s)}}$ of  the $s$-th weak derivative $f^{(s)}$ of a function $f$, if it exists, is also equivalently given by a specific weighted norm $\norm_\hw$ with an appropriate choice of weights $\hw$ (c.f. \cite{Neubauer1988}). Thus, by considering the corresponding $\cF_{\hw}$-risk we also cover the estimation of derivatives of the slope function. This question is also discussed in detail in Section \ref{sec:ex}.

In this paper we characterize the a-priori information on the slope parameter such as smoothness by considering ellipsoids (see definition below) in $L^2[0,1]$ with respect to a weighted norm $\norm_\bw$ for a pre-specified weight sequence $\bw$. Again an appropriate choice of the sequence $\bw$ enables us not only to restrict the slope parameter to a class of differentiable functions (considered, e.g. in \cite{CrambesKneipSarda2007}) but, for instance, also to a class of analytic functions. Moreover, it is usually  assumed that the  sequence $\lambda$  of eigenvalues of $\op$ has a polynomial decay (c.f. \cite{HallHorowitz2007} or \cite{CrambesKneipSarda2007}). However, it is well-known that this restriction may exclude several interesting cases, such as an exponential decay. Therefore, we do not  impose a specific  form of a decay, but  consider  a  third sequence of weights $\lw$ characterizing the decay of $\lambda$. Then we show that the three sequences $\bw$ (regularity of $\beta$), $\lw$ (regularity of $\op$) and $\hw$ (measure of the performance of the estimator) determine together the obtainable accuracy of any estimator. In other words, in Section \ref{sec:gen} we derive a lower bound under minimal regularity conditions on these sequences. It is remarkable, that  a simple orthogonal series estimator attains this lower bound up to a constant under very mild moment assumptions on the regressor and the error term. 

To be more precise, we replace the unknown quantities $g_j$ and $\lambda_j$ in equation \eqref{intro:e2} by their empirical counterparts. That is, if  $(Y_1,X_1),\dotsc,(Y_n,X_n)$ denotes an i.i.d. sample of $(Y,X)$, then for each $j\geqslant 1$, we consider the unbiased estimator
\begin{equation}\label{bm:def:est:EV-g}
\widehat{g}_j:=\frac{1}{n}\sum_{i=1}^n Y_{i}\,\skalarV{X_i,\psi_j},\quad\mbox{and}\quad \widehat{\lambda}_{j}:=\frac{1}{n}\sum_{i=1}^n \skalarV{X_{i},\psi_{j}}^2\end{equation}
for $g_j$ and $\lambda_j$ respectively. The orthogonal series estimator $\widehat{\beta}$ of $\beta$  is then defined  by
\begin{equation}\label{bm:def:est:reg}
\widehat{\beta}:=\sum_{j=1}^m \frac{\widehat{g}_{j}}{\widehat{\lambda}_{j}} \cdot\1\{\widehat{\lambda}_{j}\geqslant \alpha\}\cdot \psi_{j},\end{equation} 
where  the dimension parameter $m=m(n)$  and the threshold $\alpha=\alpha(n)$ has to tend to infinite and zero respectively as the sample size $n$  increases. Note that we introduce an additional threshold $\alpha$ on each estimated eigenvalue $\widehat{\lambda}_{j}$, since it could be arbitrarily close to zero  even in case that the true eigenvalue $\lambda_j$ is sufficiently far away from zero.  Thresholding in the Fourier domain  has been used, for example, in a deconvolution problem in  \cite{MairRuymgaart96}, \cite{Neumann1997} or \cite{Jo07} and coincides with an approach called spectral cut-off in the numerical analysis literature  (c.f. \cite{Tautenhahn96}).

The paper is organized in the following way.  In Section \ref{sec:methodology} we formalize the regularity conditions  on the slope parameter $\beta$ 
 and the covariance operator $\op$ characterized through different weight sequences. Moreover, we state the minimal conditions on these weight sequences as well as the moments of the random function $X$ and the error term $\epsilon$ used throughout the paper. In Section 3 we show consistency in the $\cF_\hw$-risk of the proposed  orthogonal series estimator under very mild assumptions. For example,  considering the $L^2$-risk, i.e., $\hw\equiv1$, there are no additional regularity conditions on the slope parameter needed. Furthermore, we derive a lower and an upper bound for the $\cF_\hw$-risk only supposing the minimal conditions on the sequences $\bw$, $\hw$ and $\lw$. These results are illustrated in Section 4 by considering the mean squared prediction error as well as the optimal estimation of derivatives of $\beta$ in case that the slope function belongs to a Sobolev space of periodic functions and  that
the covariance operator  $T_{\cov}$ is finitely or infinitely smoothing. All proofs can be found in the Appendix.

\section{Notations and basic assumptions}\label{sec:methodology}
\paragraph{Second order stationarity.}In this paper we suppose that  the regressor $X$ is second order stationary, i.e.,  there exists a positive definite function 
 $c:[-1,1]\to\R$  such that $\cov(X(t),X(s))=c(t-s)$, $s,t\in[0,1]$. Thereby we show in Proposition \ref{app:prop:cov} in the Appendix that the eigenfunctions of the covariance operator $T_{\cov}$ associated to $X$ are given by the trigonometric basis
\begin{equation}\label{bm:def:trigon}
\psi_{1}:\equiv1, \;\psi_{2j}(s):=\sqrt{2}\cos(2\pi j s),\; \psi_{2j+1}(s):=\sqrt{2}\sin(2\pi j s),s\in[0,1],\; j\in\N\end{equation}
and the corresponding eigenvalues satisfy
\begin{equation}\label{bm:def:EV}\lambda_{1}=\int_{-1}^1c(s)ds,\quad \lambda_{2j}=\lambda_{2j+1}=\int_{-1}^1\cos(2\pi j s) c(s)ds,\; j\in\N. \end{equation}
Notice that the eigenfunctions are known to the statistician and only the eigenvalues depend on the unknown covariance function $c(\cdot)$, i.e., have to be estimated. 

\paragraph{Minimal regularity conditions.}  It is well-known that the obtainable accuracy of any estimator of the slope parameter $\beta$ is essentially determined by additional regularity conditions imposed   on both  the slope parameter $\beta$ and the sequence of eigenvalues $(\lambda_j)$ of the covariance operator. In this paper these conditions are characterized through different weighted norms in $L^2[0,1]$, which we formalize  now. Given a  strictly positive sequence of  weights $w:=(w_j)_{j\geqslant1}$ and  a constant $c>0$ denote for all $r\in\R$ by $\cF_{w^r}^c$ the ellipsoid given by \begin{equation*}
 \cF_{w^r}^c := \Bigl\{f\in L^2[0,1]: \sum_{j=1}^\infty w_j^r |\skalarV{f,\psi_j}|^2=:\normV{f}_{w^r}^2\leq c\Bigr\}.
\end{equation*}
Furthermore, let $\cF_{w^r}:=\{f\in L^2[0,1]: \normV{f}_{w^r}^2<\infty\}$.  Here and subsequently, we suppose that given a strictly positive sequence of weights $\bw:=(\bw_j)_{j\geqslant1}$ the slope  function $\beta$ belongs to the ellipsoid $\cF_\bw^\br$ for some $\br>0$. The ellipsoid $\cF_\bw^\br$  captures then all the prior information (such as smoothness) about the unknown slope function $\beta$. It is worth to note, that  in case $\bw\equiv 1$ the set $\cF_\bw^\br$ denotes an ellipsoid in $L^2[0,1]$ and hence  does not imposes additional restrictions on $\beta$. Furthermore, given a  strictly positive sequence of weights $\lw:=(\lw_j)_{j\geqslant1}$    we assume that the sequence of eigenvalues $(\lambda_j)_j$ of the covariance operator $\op$ is an element of  the  set $\cS_{\lw}^\ld$  defined for $\ld\geq 1$ by
\begin{equation}\label{res:link:gen}
\cS_{\lw}^\ld:=\Bigl\{ (\lambda_j)_{j\geqslant 1}: 1/d \leqslant \lambda_j /\upsilon_j\leqslant d, \quad\forall j\in\N\Bigr\}.
\end{equation}
Notice  that the sequence of eigenvalues  $(\lambda_j)_{j\geqslant 1}$ is  summable, since $\sum_{j\in\N} \lambda_j=\Ex\normV{X}^2<\infty$.  Therefore, the sequence $\lw$ has also to be summable. We consider this quite general class 
of eigenvalues first. However, we illustrate condition \eqref{res:link:gen} in Section \ref{sec:ex} below  by assuming a \lq\lq regular decay\rq\rq\ of the eigenvalues.  Moreover, consider a strictly positive sequence of weights $\hw:=(\hw_j)_{j\geqslant1}$. Then we shall  measure the performance of an estimator $\widehat{\beta}$ of $\beta$  by the $\cF_{\hw}$-risk, that is $\Ex\|\widehat{\beta}-\beta\|_{\hw}^2$.  In Section \ref{sec:ex}  this approach is illustrated by considering different weight sequences $\hw$. Roughly speaking, an appropriate choice of  $\hw$ enables us to cover both the estimation  of derivatives of $\beta$ as well as the optimal estimation in terms of the mean prediction error. Finally, all the results below are derived under the following minimal regularity conditions.  

\begin{assumption}\label{ass:reg} Let $\hw:=(\hw_j)_{j\geqslant 1}$, $\bw:=(\bw_j)_{j\geqslant 1}$ and  $\lw:=(\lw_j)_{j\geqslant 1}$ be strictly positive sequences of weights  with $\hw_1= 1$,  $\bw_1= 1$ and  $\lw_1= 1$ such that   $\bw$ and     $(\bw_j/\hw_j)_{j\geqslant 1}$ are nondecreasing and  $\lw$ is nonincreasing  with $\Lw:= \sum_j \lw_j< \infty$. \end{assumption}
Note that under Assumption \ref{ass:reg} the ellipsoid  $\cF_\bw^\br$ is a subset of $\cF_\hw^\br$, and hence the $\cF_{\hw}$-risk a well-defined risk for $\beta$. Roughly speaking, if $\cF_\bw^\br$ describes $p$-times differentiable functions, then the Assumption \ref{ass:reg} ensures that the $\cF_{\hw}$-risk involves maximal $s\leqslant p$ derivatives.

\paragraph{Moment assumptions.}The results  derived below involve  additional  conditions on the moments of the  random function $X$ and the error term $\epsilon$, which we formalize now. Let $\cX$ be the set of all centered  second order stationary random functions $X$ with finite second moment, i.e., $\Ex\normV{X}^2<\infty$, and strictly positive covariance operator. Then given  $X\in\cX$  the random variables $\{\skalarV{X,\psi_j}/\sqrt{\lambda_j}, j\in\N\}$  are  centered with variance one and moreover pairwise  uncorrelated.  Here and subsequently,  $\cX^{m}_{\eta}$, $m\in\N$, $\eta\geq1$, denotes the subset of $\cX$  containing all random functions $X$ such that the $m$-th moment of the corresponding standardized  random variables  $\{\skalarV{X,\psi_j}/\sqrt{\lambda_j}, j\in\N\}$  are  uniformly bounded, that is
\begin{equation}\label{bm:def:X}
\cX^{m}_{\eta}:=\Bigl\{ X\in\cX\;\text{ with }\sup_{j\in\N} \Ex\Bigl|\frac{\skalarV{X,\psi_j}}{\sqrt{\lambda_j}}\Bigr|^m \leq \eta \Bigr\}.
\end{equation}
It is worth noting that in case $X\in\cX$ is a Gaussian random function the corresponding random variables  $\{\skalarV{X,\psi_j}/\sqrt{\lambda_j}, j\in\N\}$ form an  i.i.d. sample of  Gaussian random variables with mean zero and variance one. Hence, for each $k\in\N$ there exists $\eta$ such that any Gaussian random function $X\in \cX$ belongs also to $\cX^{k}_{\eta}$. In what follows, $\cE^m_\eta$ stands for the set of all centered error terms $\epsilon$ with variance one and finite $m$-th moment, i.e., $\Ex|\epsilon|^m\leq \eta$.

\section{Optimality in the general case}\label{sec:gen}
\paragraph{Consistency.}The $\cF_{\hw}$-risk of the estimator $\widehat\beta$ given in (\ref{bm:def:est:reg}) is essentially determined by the deviation of the estimators of $(g_j)_j$ and $(\lambda_j)_j$  and  by the regularization error due to the threshold.
The next assertion summarizes  minimal conditions to ensure consistency of the estimator defined in (\ref{bm:def:est:reg}). 
 \begin{prop}[Consistency]\label{res:prop:cons} Assume an $n$-sample of $(Y,X)$ satisfying \eqref{intro:e1} with $\sigma>0$. Let $\beta\in \cF_\bw$, $X\in\cX^4_{\eta}$ and $\epsilon\in\cE^4_\eta$, $\eta\geqslant 1$.  Consider the estimator $\widehat{\beta}$
 with threshold  $m:=m(n)$ and parameter $\alpha:=\alpha(n)$ satisfying $m\to\infty$, $\alpha=o(1)$ and   $(\sup_{j \leqslant m} \hw_j) (n\alpha^2)^{-1}=o(1)$   as  $n\to\infty$. If in addition $\bw$ and $\hw$ satisfy Assumption \ref{ass:reg}, then
$\Ex\|\widehat{\beta}-\beta\|_{\hw}^2=o(1)$ as  $n\to\infty$. 
 \end{prop}
\begin{rem}Since the last result covers the case $\bw\equiv\hw \equiv  1$ it follows that  the estimator  $\widehat\beta$ is consistent without any additional restriction on $\beta\in L^2[0,1]$ provided $m\to\infty$, $\alpha=o(1)$ and   $n\alpha^2\to\infty$   as  $n\to\infty$.
\hfill$\square$\end{rem}

\paragraph{The lower bound.}
It is well-known  that in general the hardest one-dimensional subproblem does not capture the full difficulty in estimating the solution of an inverse problem even in case of a known operator (for details see e.g. the proof in \cite{MairRuymgaart96}).  In other words, there does not exist two sequences of slope functions $\beta_{1,n},\beta_{2,n}\in \cF_\bw^\br$, which are statistically not consistently distinguishable and satisfy $\normV{\beta_{1,n}-\beta_{2,n}}^2_\hw\geqslant C \dstar$, where $\dstar$ is the optimal rate of convergence. Therefore we need to consider subsets of $\cF_\bw^\br$ with  growing number of elements in order to get the optimal lower bound. More specific, we obtain the following lower bound by applying Assouad's cube technique (see e.g. \cite{KorostelevTsybakov1993} or \cite{ChenReiss2008})   under the additional assumption that the error term $\epsilon$ is standard normal distributed, i.e., $\epsilon\sim\cN(0,1)$, and independent of the regressor. 
\begin{theo}\label{res:lower}Assume an $n$-sample of $(Y,X)$ obeying \eqref{intro:e1} with $\sigma>0$.    Suppose that the error term $\epsilon\sim\cN(0,1)$ is independent of the second order stationary regressor $X$ with associated sequence of eigenvalues  $(\lambda_j)\in \cS_\lw^\ld$.  Consider $\cF_\bw^\br,$ $\br>0,$  as set of slope functions. Let $\kstar:=\kstar(n)\in\N$ and $\dstar:=\dstar(\kstar)\in\R^+$ for some $\triangle\geq 1$  be chosen such that
\begin{equation}\label{res:lower:def:md}
1/\triangle\leq \frac{\bw_{\kstar}}{n\,\hw_{\kstar}}\sum_{j=1}^{\kstar} \frac{\hw_{j}}{\lw_j} \leq \triangle\quad\text{ and }\quad \dstar:=\hw_{\kstar}/\bw_{\kstar}.
\end{equation}
If in addition the Assumption \ref{ass:reg} is satisfied then for any estimator $\breve\beta$ we have 
\begin{equation*}  \sup_{\beta \in \cF_\bw^\br} \left\{ \Ex\normV{\breve{\beta}-\beta}^2_\hw\right\}\geqslant  \frac{1}{4\,\triangle}\min \left( \frac{\sigma^2}{2\,d},  \frac{\br}{\triangle}\right) \,\max(\dstar,1/n).
\end{equation*}
\end{theo}
\begin{rem}\label{rem:lower}The normality assumption in the last theorem is only used to simplify the calculation of the distance  between distributions corresponding to different slope functions. Obviously the derived lower bound is  still valid if we consider the less restrictive assumption  that  the error term $\epsilon$ belongs to $\cE_\eta^m$ for some $m\in\N$ and sufficiently large $\eta$. Furthermore, it is worth to note that the lower bound tends only to zero if  $(\hw_{j}/\bw_{j})$ is a zero sequence. In other words, in case $\bw\equiv 1$, i.e., without  any additional restriction on $\beta\in L^2[0,1]$, uniform consistency over $L^2[0,1]$ in the $\cF_\hw$-risk is only possible if the weighted norm $\norm_\hw$ is weaker than the usual $L^2$-norm, that is, $\hw$ is a zero sequence. This obviously reflects the ill-posedness of the underlying inverse problem.\hfill$\square$\end{rem}
\paragraph{The upper bound.}
The next theorem states that the rate $\max(\dstar,1/n)$ of the lower bound given  in Theorem \ref{res:lower} provides also an upper bound of 
the proposed estimator $\widehat{\beta}$. Therefore the rate $\max(\dstar,1/n)$ is optimal and hence the estimator $\widehat{\beta}$   is  minimax-optimal. 
\begin{theo}\label{res:upper}Assume an $n$-sample of $(Y,X)$ satisfying \eqref{intro:e1} with $\sigma>0$.   Suppose that the regressor $X$ is second order stationary with associated sequence of eigenvalues  $(\lambda_j)\in \cS_\lw^\ld$.  Consider $\kstar:=\kstar(n)$ and  $\dstar:=\dstar(n)$  given in \eqref{res:lower:def:md} for some $\triangle\geqslant 1$. Let  $\widehat{\beta}$  be  the estimator defined in \eqref{bm:def:est:reg}  with  $m:=\kstar$ and  $\alpha:=(1/n)\min(1,\bw_{\kstar}/(2\ld\triangle))$.
 If in  addition  $X\in\cX^{4k}_{\eta}$ and $\cE^{4k}_\eta$, $k\geqslant4$, then  for some generic constant $C>0$ we have
\begin{gather*}
\sup_{\beta \in \cF_\bw^\br}\left\{\Ex\normV{\widehat{\beta}- \beta}^2_\hw\right\}\leqslant C \, \ld^5\,\triangle^3\,\eta\,[ \br\, \ld\,  \Lw+\sigma^2]\, \max( \dstar,1/n),
\end{gather*}
for all sequences  $\bw$, $\hw$ and $\lw$ satisfying Assumption \ref{ass:reg}.
\end{theo}
\begin{rem} It is worth to note that the bound derived in the last theorem is non asymptotic. Furthermore,  as in case of the lower bound (see Remark \ref{rem:lower}) also the upper bound tends only to zero,  if $ (\hw_j/\bw_j)$ is a zero sequence.  Therefore  the estimator $\widehat{\beta}$ is consistent even without any additional restriction on $\beta\in L^2[0,1]$, i.e., $\bw\equiv 1$, as long as  $\hw$ is a zero sequence. We shall stress that from  Theorem \ref{res:upper} follows that for all sequences  $\bw$, $\hw$ and $\lw$ satisfying the minimal  regularity Assumption \ref{ass:reg} the orthogonal series estimator $\widehat{\beta}$ attains the optimal rate $\max(\dstar,1/n)$  and hence is minimax-optimal. In particular, it is easily seen that the optimal rate $\max(\dstar,1/n)$ is parametric if and only if $\sum_{j=1}^\infty \hw_j/\lw_j <\infty$. Hence, in this case the rate of the orthogonal series estimator $\widehat{\beta}$ is parametric again without any additional restriction on $\beta\in L^2[0,1]$, i.e., $\bw\equiv 1$. Finally as long as the sequence $\bw$ is unbounded in Theorem \ref{res:upper} the threshold parameter $\alpha$ satisfies $\alpha=1/n$ for all  sufficiently large $n$. Thus in this situation as open problem remains only how to choose the dimension parameter $m$  adaptively from the data. We are currently exploring this issue.\hfill$\square$\end{rem}
\section{Mean prediction error and derivative estimation}\label{sec:ex}
In this section we suppose the slope function $\beta$ is an element of the Sobolev space of periodic functions $\cW_p$  
 given for $p> 0$ by  \begin{equation*}
 \cW_{p}=\Bigl\{f\in H_{s}: f^{(j)}(0)=f^{(j)}(1),\quad j=0,1,\dotsc,p-1\Bigr\},
 \end{equation*}
where  $H_{p}:= \{ f\in L^2[0,1]:  f^{(p-1)}\mbox{ absolutely continuous }, f^{(p)}\in L^2[0,1]\}$  is a Sobolev space (c.f. \cite{Neubauer1988,Neubauer88}, \cite{MairRuymgaart96} or \cite{Tsybakov04}). However, if we consider the sequence of weights  $(w_j^p)_{j\in\N}$ given by
\begin{equation}\label{form:def:b}
w_{1}^p=1\quad\mbox{ and }\quad w_{2j}^p= w_{2j+1}^p=j^{2p},\qquad j\in\N.
\end{equation}
Then the Sobolev space $\cW_p$ of periodic functions  is equivalently given by $\cF_{w^p}$. Therefore, let us denote by $\cW_p^\br:=\cF_{w^p}^\br$, $\br>0$, an ellipsoid in the Sobolev space $\cW_p$. We use  in case $p=0$ again the convention that $\cW_p^\br$ denotes an ellipsoid in $L^2[0,1]$.

\paragraph{Mean prediction error.} We shall first measure the performance of an estimator $\widehat{\beta}$ by  the mean prediction error (MPE), i.e., $\Ex\skalarV{\op(\widehat\beta-\beta),(\widehat\beta-\beta)}$. Consequently, if the sequence of eigenvalues $(\lambda_j)$ associated to the covariance operator  $\op$  satisfies a link condition, that is $(\lambda_j) \in \cS_{\lw}^{\ld}$ for some weight sequence  $\lw$ (see definition \eqref{res:link:gen}). Then the MPE  is equivalent to the $\cF_\hw$-risk with $\hw\equiv\lw$, that is $\Ex \normV{\widehat\beta-\beta}_\lw^2 \asymp_d \Ex\skalarV{\op(\widehat\beta-\beta),\widehat\beta-\beta}$. To illustrate the previous results we assume in the following the sequence $\lw$ to be either  polynomially decreasing, i.e., $\lw_1=1$ and $\lw_j  = |j|^{-2a}$, $j\geqslant 2$, for some $a>1/2$,  or   exponentially decreasing, i.e.,  $\lw_1=1$ and $\lw_j  = \exp(-|j|^{2a})$, $j\geqslant 2$, for some $a>0$. In the polynomial case easy calculus shows  that  a covariance operator ${\op}$ with eigenvalues $(\lambda_j)\in\cS_\lw^\ld$, i.e., $\lambda_j\asymp_d |j|^{-2a}$,  acts like integrating   $(2a)$-times and hence it is called {\it finitely smoothing} (c.f. \cite{Natterer84}). This is the case considered, for example, in  \cite{CrambesKneipSarda2007}. On the other hand in the exponential case  it can  easily be seen that the link condition $(\lambda_j)\in\cS_\lw^\ld$, i.e., $\lambda_j\asymp_d \exp(-j^{2a})$, implies $\cR({\op})\subset \cW_{s}$ for all $s>0$, therefore  the operator ${\op}$ is called {\it infinitely smoothing} (c.f. \cite{Mair94}). Since in both cases the minimal regularity conditions given in Assumption \ref{ass:reg} are satisfied, the lower bounds presented in the next assertion follow directly from Theorem \ref{res:lower}. Here and subsequently, we write $a_n\lesssim b_n$  when there exists $C>0$ such that  $a_n\leqslant C\, b_n$  for all sufficiently large $ n\in\N$ and  $a_n\sim b_n$ when $a_n\lesssim b_n$ and $b_n\lesssim a_n$ simultaneously. 

\begin{prop}\label{MPE:lower}\dr \mbox{Under the assumptions of Theorem \ref{res:lower} we have for any estimator $ \breve{\beta}$}\\[-4ex]
\begin{itemize}\item[(i)] in the polynomial case, i.e. $\lw_1=1$ and $\lw_j  = |j|^{-2a}$, $j\geqslant 2$, for some $a>1/2$,    that\\[1ex] 
\hspace*{5ex}$\sup_{\beta \in \cW_p^\rho} \Bigl\{ \Ex\skalarVl{\op(\breve\beta-\beta),(\breve\beta-\beta)}\Bigr\}\gtrsim n^{-(2p+2a)/(2p+2a+1)} $,
\item[(ii)] in the exponential case, i.e. $\lw_1=1$ and $\lw_j  = \exp(-|j|^{2a})$, $j\geqslant 2$, for some $a>0$,    that\\[1ex] 
\hspace*{5ex}$\sup_{\beta \in \cW_p^\rho} \Bigl\{ \Ex\skalarVl{\op(\breve\beta-\beta),(\breve\beta-\beta)}\Bigr\}\gtrsim n^{-1}(\log n)^{1/2a}$.
\end{itemize}
\end{prop}

On the other hand, if the dimension parameter  $m$ and the threshold $\alpha$ in the definition of  the estimator $\widehat{\beta}$ given in \eqref{bm:def:est:reg} are chosen appropriate, then by applying Theorem \ref{res:upper}  the rates of the lower bound given in the last assertion provide up to a constant also the upper bound of the risk of the  estimator $\widehat{\beta}$, which is summarized in the next proposition. 
We have thus proved that these rates are  optimal and the proposed estimator $\widehat{\beta}$ is minimax optimal in both cases. 

\begin{prop}\label{MPE:upper}\dr Under the assumptions of Theorem \ref{res:upper} consider the estimator $\widehat{\beta}$\\[-4ex]
\begin{itemize}\item[(i)] in the polynomial case, i.e. $\lw_1=1$ and $\lw_j  = |j|^{-2a}$, $j\geqslant 2$, for some $a>1/2$,  with dimension $m\sim n^{1/(2p+2a+1)}$ and threshold $\alpha\sim 1/n$. Then we have\\[1ex] 
\hspace*{5ex}$\sup_{\beta \in \cW_p^\rho} \Bigl\{ \Ex\skalarVl{\op(\widehat\beta-\beta),(\widehat\beta-\beta)}\Bigr\}\lesssim n^{-(2p+2a)/(2p+2a+1)}$,
\item[(ii)] in the exponential case, i.e. $\lw_1=1$ and $\lw_j  = \exp(-|j|^{2a})$, $j\geqslant 2$, for some $a>0$,  with dimension  $m\sim (\log n)^{1/(2a)}$ and threshold $\alpha\sim 1/n$. Then\\[1ex] 
\hspace*{5ex}$\sup_{\beta \in \cW_p^\rho} \Bigl\{ \Ex\skalarVl{\op(\widehat\beta-\beta),(\widehat\beta-\beta)}\Bigr\}\lesssim n^{-1}(\log n)^{1/2a}$.
\end{itemize}
 \end{prop}
\begin{rem}It is of interest to compare our results with those of \cite{CrambesKneipSarda2007} who measure the performance of their estimator in terms of the prediction error. In their notations the decrease of the eigenvalues of ${\op}$ is assumed to be of order $ (|j|^{-2q-1})$, i.e., $q=a-1/2$.  Furthermore they  suppose the slope function to be $m$-times continuously differentiable, i.e., $m=p$. By using this reparametrization we see that our results in the polynomial case imply the same rate of convergence in probability of the prediction error as it is presented in  \cite{CrambesKneipSarda2007}.  However, from our general results follows a lower and an upper bound of the MPE not only in the polynomial case but also  in the exponential case. 

Furthermore, we shall emphasize the interesting influence of the parameters $p$ and $a$ characterizing the smoothness of $\beta$ and the decay of the eigenvalues of  $\op$, respectively. As we see from Propositions \ref{MPE:lower} and \ref{MPE:upper},  in the polynomial case  an increasing  value of $p$ leads to a faster optimal rate. In other words, as expected, a smoother  regression function can be faster estimated. 
The situation in the exponential case is extremely different. It seems rather surprising that, contrary  to the polynomial case, in the exponential case the optimal rate of convergence does not depend on the value of $p$, however this dependence is clearly hidden in the constant. Furthermore, the dimension parameter $m$  does not even depend on the value of $p$.  Thereby, the proposed estimator is automatically adaptive, i.e., it does not involve an a-priori knowledge of the degree of smoothness of the slope function $\beta$. However, the choice of the dimension parameter depends on  the value $a$ specifying the decay of the eigenvalues of  $\op$. Note further that in  both cases 
 an increasing value of $a$ leads to a faster optimal rate of convergence, i.e., we may call $1/a$  {\it degree of ill-posedness} (c.f. \cite{Natterer84}).   Finally, we shall stress that Proposition \ref{MPE:upper} covers the case $p=0$, i.e., $\widehat{\beta}$ is consistent with optimal MPE-rate without  additional restrictions on $\beta\in L^2[0,1]$.$\square$\end{rem}
\paragraph{Estimation of the derivatives.} 
Let us consider now the estimation of derivatives of the slope function $\beta$. It is well-known, that for any function $g$ belonging to a  Sobolev-ellipsoid $\cW_{p}^\br=\cF_{w^p}^\br$ with weights $w^p$ given in \eqref{form:def:b} the weighted norm  $\normV{g}_{w^s}$ for each $0\leqslant s\leqslant p$  is equivalent to the $L^2$-norm of the $s$-th weak derivative $g^{(s)}$, that is,   $\normV{g^{(s)}}\asymp_{(2\pi)^{2s}}\normV{g}_{w^s}$. Thereby, the results in the Section \ref{sec:gen} imply again a lower bound as well as an upper bound of the $L^2$-risk for the estimation of the  $s$-th weak derivative of  $\beta$. In the following we consider again the two particular cases of polynomial and exponential decreasing rates for the sequence of weights $(\lw_j)$. The next assertion summarizes then lower bounds for the $L^2$-risk for the  estimation  of the $s$-th weak derivative $\beta^{(s)}$ of $\beta$ in both cases. 

\begin{prop}\label{MSE:lower}\dr \mbox{Under the assumptions of Theorem \ref{res:lower} we have for any estimator $ \widetilde{\beta}^{(s)}$}\\[-4ex]
\begin{itemize}\item[(i)] in the polynomial case, i.e. $\lw_1=1$ and $\lw_j  = |j|^{-2a}$, $j\geqslant 2$, for some $a>1/2$,    that\\[1ex] 
\hspace*{5ex}$\sup_{\beta \in \cW_p^\rho} \bigl\{ \Ex\normV{\widetilde{\beta}^{(s)}-\beta^{(s)}}^2\bigr\}\gtrsim n^{-(2p-2s)/(2p+2a+1)} $,
\item[(ii)] in the exponential case, i.e. $\lw_1=1$ and $\lw_j  = \exp(-|j|^{2a})$, $j\geqslant 2$, for some $a>0$,    that\\[1ex] 
\hspace*{5ex}$\sup_{\beta \in \cW_p^\rho} \bigl\{ \Ex\normV{\widetilde{\beta}^{(s)}-\beta^{(s)}}^2\bigr\}\gtrsim (\log n)^{-(p-s)/a}$.
\end{itemize}
\end{prop}

On the other hand considering the estimator $\widehat\beta$ given in \eqref{bm:def:est:reg}, we only have to calculate the $s$-th derivative of $\widehat\beta$. However, given the exponential basis, which is linked to the trigonometric basis by the relation $\exp(2 \iota \pi k  t) = 2^{-1/2} ( \psi_{2k}(t) + \iota \ \psi_{2k+1}(t)),$ for  $k \in \Z$ and $t \in [0,1],$ with $\iota^2=-1,$ then   for $0\leqslant s <p$ the $s$-th derivative $\widehat\beta^{(s)}$ of $\widehat\beta$ in a weak sense  is
\begin{equation}\label{bm:def:est:reg:s}
\widehat\beta^{(s)}(t) =\sum_{k\in\Z} ( 2\iota \pi k)^s  \left( \int_0^1 \widehat\beta(u) \exp(-2 \iota \pi k u) \ du \right)  \exp(2 \iota  \pi k t),\quad t \in [0,1].
\end{equation}
Note, that the sum in \eqref{bm:def:est:reg:s} contains only a finite number of nonzero summands and hence its numerical implementation is straightforward. 
Furthermore,  if the dimension parameter  $m$ and the threshold $\alpha$ in the definition of  the estimator $\widehat{\beta}$ given in \eqref{bm:def:est:reg} are chosen appropriate, then by applying Theorem \ref{res:upper}  the rates of the lower bound given in the last assertion provide up to a constant again the upper bound of the $L^2$-risk of the  estimator $\widehat{\beta}^{(s)}$, which is summarized in the next proposition. We have thus proved that these rates are  optimal and the proposed estimator $\widehat{\beta}^{(s)}$ is minimax optimal in both cases.
\begin{prop}\label{MSE:upper}\dr Under the assumptions of Theorem \ref{res:upper} consider the estimator $\widehat{\beta}^{(s)}$\\[-4ex]
\begin{itemize}\item[(i)] in the polynomial case, i.e. $\lw_1=1$ and $\lw_j  = |j|^{-2a}$, $j\geqslant 2$, for some $a>1/2$,  with $m\sim n^{1/(2p+2a+1)}$ and threshold $\alpha\sim n$. Then\\[1ex] 
\hspace*{5ex}$\sup_{\beta \in \cW_p^\rho} \bigl\{ \Ex\normV{\widehat{\beta}^{(s)}-\beta^{(s)}}^2\bigr\}\lesssim n^{-(2p-2s)/(2p+2a+1)}$,
\item[(ii)] in the exponential case, i.e. $\lw_1=1$ and $\lw_j  = \exp(-|j|^{2a})$, $j\geqslant 2$, for some $a>0$,  with $m\sim (\log n)^{1/(2a)}$ and threshold $\alpha\sim n$. Then\\[1ex] 
\hspace*{5ex}$\sup_{\beta \in \cW_p^\rho} \bigl\{ \Ex\normV{\widehat{\beta}^{(s)}-\beta^{(s)}}^2\bigr\}\lesssim (\log n)^{-(p-s)/a}$.
\end{itemize}
 \end{prop}

\begin{rem}\label{rem:sob:upper:pol:1}It is worth noting that the $L^2$-risk in estimating the slope function $\beta$ itself, i.e., $s=0$, has been considered in \cite{HallHorowitz2007} only in the polynomial case. In their notations the decrease of the eigenvalues of ${\op}$ is of order $ (|j|^{-\alpha})$, i.e., $\alpha=2a$.  Furthermore the Fourier coefficients of  the slope function decay at least with rate $j^{-\beta}$, i.e., $\beta=p+1/2$. By using this reparametrization we see that we recover the result of \cite{HallHorowitz2007} in the polynomial case with $s=0$, but without the additional assumption $ \beta > \alpha/2+1$ or $\beta >\alpha-1/2$. 

Furthermore, we shall discuss again  the influence of the parameters $p$, $s$ and $a$. As we see from Propositions \ref{MSE:lower} and \ref{MSE:upper}, in both cases an  decreasing of the value of $a$ or an increasing of the value $p$  leads to a faster  optimal rate of convergence.
Hence, in opposite to the MPE  by considering the $L^2$-risk the parameter $a$ describes in both cases the {\it degree of ill-posedness}. Furthermore, the estimation of higher derivatives of the slope function, i.e. by considering a larger value of $s$, is as usual only possible with a slower optimal rate. Finally,  as for the MPE  in the exponential case the dimension parameter $m$  does not depend on the values of $p$ or $s$, hence the proposed estimator is automatically adaptive. \hfill$\square$\end{rem}

\begin{rem}\label{rem:sob:upper:pol:2}There is an interesting issue hidden in the parametrization we have chosen. Consider a classical indirect regression model with known operator given by ${\op}$, i.e., $Y=[{\op}\beta](U)+\epsilon$ where $U$ has a uniform distribution on $[0,1]$ and $\epsilon$ is white noise  (for details see e.g. \cite{MairRuymgaart96}). If in addition the operator ${\op}$ is finitely smoothing, i.e., $(\lw_j)$ is polynomially decreasing with $\lw_j=j^{-2a}$, $j\geqslant 2$. Then given an $n$-sample of $Y$ the optimal rate of convergence of the $L^2$-risk of any estimator of $\beta^{(s)}$ is of order $n^{-2(p-s)/[2(p+2a)+1]}$,  since $\cR({\op})=\cW_{2a}$ (c.f. \cite{MairRuymgaart96} or \cite{ChenReiss2008}). However, we have shown that in a functional linear model even with estimated operator the optimal rate is of order $n^{-2(p-s)/[2(p+a)+1]}$. Thus comparing both rates we see that in a functional linear model the covariance operator ${\op}$ has  the {\it degree of ill-posedness} $a$ while the same operator has in the indirect regression model a {\it degree of ill-posedness} $(2a)$. In other words in a functional linear model we do not face the complexity of an inversion of ${\op}$ but only of its square root ${\op}^{1/2}$. This, roughly speaking, may be seen as a multiplication of the normal equation $YX=\skalarV{\beta,X}X+X\epsilon$ by the inverse of $\op^{1/2}$. Notice that ${\op}$ is also the covariance operator associated to the error term  $\epsilon X$. Thus the multiplication by the inverse of $\op^{1/2}$ leads, roughly speaking, to white noise  and hence  to an indirect regression model rather defined  by $T_{\cov}^{1/2}$ than $T_{\cov}$. The same finding holds true in case of an infinitely smoothing operator ${\op}$. However, in this situation $(\log n)^{-(p-s)/a}$ is the optimal rate  in an indirect regression model  given by ${\op}$ as well as  $\op^{1/2}$. Thus, the above described effect is not visible formally, but is actually hidden in the order symbol.\hfill$\square$\end{rem}
\appendix
\section{Appendix}\label{app:proofs}
\begin{prop}\label{app:prop:cov}Let $X$ be second order stationary with $\Ex [X(t)X(s)]=c(t-s)$, $t,s\in[0,1]$, for some positive definite function $c:[-1,1]\to\R$. Then the associated covariance operator $T_{\cov}$ admits an eigenvalue decomposition with eigenfunctions given by the trigonometric basis defined in \eqref{bm:def:trigon} and  corresponding eigenvalues given by \eqref{bm:def:EV}.
 \end{prop}
\noindent\textcolor{darkred}{\sc Proof.}
Let $f\in L^2[0,1]$ and consider  $g=T_{\cov}f=\int_{0}^1 f(t) c(\cdot-t)dt$. Since $c$ is even, it is straightforward to show that $\int_{0}^1g(s)e^{-is\lambda}ds=\int_{0}^1f(s)e^{-is\lambda}ds\int_{-1}^1c(s)\cos(s\lambda)ds$ and  $\int_{0}^1g(s)e^{is\lambda}ds=\int_{0}^1f(s)e^{is\lambda}ds\int_{-1}^1c(s)\cos(s\lambda)ds$ for all $\lambda\in\R$. Due to  this  we obtain  for all $\lambda\in\R$ the following identities 
\begin{gather*} \int_{0}^1g(s)\cos(s\lambda)ds=\int_{0}^1f(s)\cos(s\lambda)ds\int_{-1}^1c(s)\cos(s\lambda)ds,\\
\int_{0}^1g(s)\sin(s\lambda)ds=\int_{0}^1f(s)\sin(s\lambda)ds\int_{-1}^1c(s)\cos(s\lambda)ds.\end{gather*}
Consider the trigonometric basis $\{\psi_{n}\}$ and the values $\{\lambda_{n}\}$ given in \eqref{bm:def:trigon} and \eqref{bm:def:EV}, respectively, then we have just shown, that
$\skalarV{T_{\cov}f,\psi_{n}}=\skalarV{f,\psi_{n}} \lambda_{n} $ for all $f\in L^2[0,1]$ and $n\in\N$, which proves the result.\hfill$\square$

\subsection{Proofs of Section \ref{sec:gen}}\label{app:proofs:gen}

We begin by defining and recalling notations to be used in the proofs:
\begin{multline}\label{app:l:def}
X_{ij}:= \skalarV{X_i,\psi_j}, \quad \beta_j=\skalarV{\beta,\psi_j},\quad
T_{n,j}:=\frac{1}{n}\sum_{i=1}^n (Y_iX_{ij} -  X_{ij}^2 \beta_j),\quad \lambda_j=\Ex {X}_{ij}^2,\\
\widetilde{\beta}_m:=\sum_{j=1}^m \beta_{j} \cdot\1\{\widehat{\lambda}_{j}  \geqslant \alpha\}\cdot \psi_{j},\quad{\beta}_m:=\sum_{j=1}^m \beta_{j} \cdot \psi_{j}.\hfill\end{multline} 
We shall prove in the end of this section two technical Lemma (\ref{app:l1:upp} - \ref{app:l1:lower}) which are used in the following proofs.

\paragraph{Proof of consistency.}\hfill\\[1ex]
\noindent\textcolor{darkred}{\sc Proof of Proposition \ref{res:prop:cons}.}
The proof is based on the  decomposition \begin{equation}\label{res:dec}
 \Ex\|\widehat{\beta}-\beta\|_{\hw}^2\leq 2\{\Ex\normV{\widehat{\beta}-
 \widetilde{\beta}_m}^2_{\hw}+\Ex\normV{\widetilde{\beta}_m- \beta}^2_{\hw}\}.
\end{equation}
We show  below under the moment condition $X\in\cX_{\eta}^4$ defined in \eqref{bm:def:X} and $\epsilon\in\cE^4_\eta$ for some universal constant $C>0$ the following  bound
\begin{equation}\label{pr:res:prop:cons:e1}
\Ex\normV{\widehat{\beta}- \widetilde{\beta}_m}^2_\hw\leqslant  C \,(\sup_{j \leqslant m} \hw_j) \,(n\alpha^2)^{-1}\, \Ex\normV{X}^{2 } \,\{\sigma^2+\normV{\beta}^2\,\Ex\normV{X}^{2}\}\,\eta,
\end{equation}
while given $\normV{\beta}_\hw<\infty$  we conclude from Lebesgue's dominated convergence theorem%
\begin{equation}\label{pr:res:prop:cons:e2}
\Ex\normV{\widetilde{\beta}_m-\beta}^2_\hw=o(1)\mbox{ in case that  }1/m=o(1),\, \alpha=o(1) \mbox{ as } n\to\infty.
\end{equation}
Thereby, the
conditions on  $m$  and $\alpha$ ensure the convergence to zero of  the two  terms on the right hand
side in \eqref{res:dec} as $n\to\infty$, which gives the result.

Proof of (\ref{pr:res:prop:cons:e1}). By  making use of the notations given in \eqref{app:l:def} it follows that
\begin{equation*}\Ex\normV{\widehat{\beta}- \widetilde{\beta}_m}^2_{\hw}=\sum_{j=1}^m  \hw_j\,\Ex \frac{|\widehat{g}_j-\beta_j\,\widehat{\lambda}_j|^2}{\widehat{\lambda}_j^2}\1\{\widehat{\lambda}_j\geqslant \alpha\}\leq \frac{1}{\alpha^2} \sum_{j=1}^m  \hw_j \Ex |T_{n,j}|^2
\end{equation*}
 and hence by using \eqref{app:l1:upp:e1}  in Lemma \ref{app:l1:upp}  we obtain  (\ref{pr:res:prop:cons:e1}).

The proof of (\ref{pr:res:prop:cons:e2}) is based on the  decomposition 
\begin{equation*}\Ex\normV{\widetilde{\beta}_m- \beta}^2_{\hw}\leqslant 2\Bigl\{ \sum_{j=1}^\infty \hw_j \, \beta_j^2\1\{j> m\} + \sum_{j=1}^m \hw_j \beta_j^2 \, P( \widehat{\lambda}_j<\alpha)\Bigr\}\leqslant 2\sum_{j=1}^\infty \hw_j \,  \beta_j^2=\normV{\beta}_\hw^2<\infty.
\end{equation*}
Thus Lebesgue's dominated convergence theorem  implies the result since in case $1/m=o(1)$ and $\alpha=o(1)$ as $n\to\infty$ for each $j\in\N$ $\1\{j> m\}=0$ and  $P( \widehat{\lambda}_j<\alpha)=o(1)$, which can be realized as follows. 
By using that $\alpha=o(1)$ as $n\to\infty$ there  exists $n_j>0$ such that for all $n \geqslant n_j$ it holds $\lambda_j\geqslant 2 \alpha$ and hence $P( \widehat{\lambda}_j<\alpha)\leqslant  P(\widehat{\lambda}_j/\lambda_j <1/2) $ together with \eqref{app:l1:upp:e3} in Lemma \ref{app:l1:upp} implies the assertion, which completes the proof. \hfill$\square$\\
 
\paragraph{Proof of the lower bound.}\hfill\\[1ex]
\noindent\textcolor{darkred}{\sc Proof of Theorem \ref{res:lower}.}  Let $X_i$, $i\in\N$, be i.i.d. copies of $X$ which is second order stationary with associated sequence of eigenvalues $(\lambda_j)_{j\geqslant 1} \in \cS_{\lw}^\ld$.   Consider independent error terms $\epsilon_i\sim \cN(0,1)$, $i\in\N$, which are  independent of the random functions $\{X_i\}$. Let  $\theta\in\{-1,1\}^{\kstar}$, where $\kstar:=\kstar(n)\in\N$ satisfies \eqref{res:lower:def:md} for some $\triangle\geqslant 1$. Consider the  $\kstar$-vector $b$ of coefficients $b_{j}$  given in  \eqref{app:l1:lower:b} in Lemma \ref{app:l1:lower}. For each $\theta$ define a slope function $\beta_\theta:=\sum_{j=1}^{\kstar}\theta_{j} b_{j}\psi_j$ which belongs to $\cF_{\bw}^\br$ due to \eqref{app:l1:lower:e1}  in Lemma \ref{app:l1:lower}. Consequently, for each $\theta$ the random variables $(Y_i,X_i)$ with $Y_i:=\int_0^1\beta_{\theta}(s)X_i(s)ds+\sigma\epsilon_i$, $i=1,\dotsc,n$, form a sample of the model \eqref{intro:e1} and we denote its joint distribution by  $P_{\theta}$. Furthermore, for $j=1,\dotsc,\kstar$ and each $\theta$ we introduce $\theta^{(j)}$ by $\theta^{(j)}_{l}=\theta_{l}$ for $j\ne l$ and $\theta^{(j)}_{j}=-\theta_{j}$.  As in case of  $P_\theta$ the conditional distribution of $Y_i$ given $X_i$  is Gaussian with mean  $\sum_{j=1}^{\kstar} \theta_{j} b_{j}X_{ij}$  and variance $\sigma^2$ it is easily seen that   the log-likelihood of $P_{\theta^{(j)}}$ w.r.t. $P_{\theta}$ is given by 
\begin{equation*}
\log\Bigl(\frac{dP_{\theta^{(j)}}}{dP_{\theta}}\Bigr)=-\frac{1}{\sigma^2}   \sum_{i=1}^n 
\Bigl\{Y_i - \sum_{l=1}^{\kstar} \theta_{l} {b}_{l} X_{il} \Bigr\} \theta_{j} {b}_{j} X_{ij} - \frac{2}{\sigma^2}   \sum_{i=1}^n b_{j}^2 X_{ij}^2
\end{equation*}
and its expectation w.r.t. $P_{\theta}$ satisfies 
$\Ex_{P_{\theta}}[\log(dP_{\theta^{(j)}}/dP_{\theta})]= -(2n/ \sigma^2) \, b^2_{j}\, \Ex X_{1j}^2 $. 
In terms of  Kullback-Leibler divergence this means $KL(P_{\theta^{(j)}},P_{\theta})= (2n/ \sigma^2) \, b^2_{j}\, \Ex X_{1j}^2 \leqslant  (2\,d\,n/ \sigma^2)\, b^2_{j}\,\lw_j$  by using that $(\lambda_j)_{j\geqslant 1} \in \cS_{\lw}^\ld$. Since the
 Hellinger distance $H(P_{\theta^{(j)}},P_{\theta})$ satisfies $H^2(P_{\theta^{(j)}},P_{\theta}) \leqslant KL(P_{\theta^{(j)}},P_{\theta})$  it follows from  \eqref{app:l1:lower:e1} in Lemma \ref{app:l1:lower} that 
\begin{equation}\label{pr:lower:e3}
H^2(P_{\theta^{(j)}},P_{\theta}) \leqslant \frac{2\,d\, n }{\sigma^2}\cdot b^2_{j}\cdot \lw_j\leqslant 1,\quad j=1,\dotsc,\kstar.
\end{equation} 
Consider  the  Hellinger affinity $\rho(P_{\theta^{(j)}},P_{\theta})= \int \sqrt{dP_{\theta^{(j)}}dP_{\theta}}$, then  for any estimator $\widetilde\beta$ follows%
\begin{align}\nonumber
\rho(P_{\theta^{(j)}},P_{\theta})&\leqslant \int \frac{|\skalarV{\widetilde{\beta}-\beta_{\theta^{(j)}},\psi_j}|}{|\skalarV{\beta_{\theta}-\beta_{\theta^{(j)}},\psi_j}|} \sqrt{dP_{\theta^{(j)}}dP_{\theta}} + \int \frac{|\skalarV{\widetilde{\beta}-\beta_{\theta},\psi_j}|}{|\skalarV{\beta_{\theta}-\beta_{\theta^{(j)}},\psi_j}|} \sqrt{ dP_{\theta^{(j)}}dP_{\theta}}\\\label{pr:lower:e4}
&\leqslant \Bigl( \int  \frac{|\skalarV{\widetilde{\beta}-\beta_{\theta^{(j)}},\psi_j}|^2}{|\skalarV{\beta_{\theta}-\beta_{\theta^{(j)}},\psi_j}|^2} dP_{\theta^{(j)}}\Bigr)^{1/2} +   \Bigl( \int  \frac{|\skalarV{\widetilde{\beta}-\beta_{\theta},\psi_j}|^2}{|\skalarV{\beta_{\theta}-\beta_{\theta^{(j)}},\psi_j}|^2} dP_{\theta}\Bigr)^{1/2}.
\end{align}
Due to the identity $\rho(P_{\theta^{(j)}},P_{\theta})=1-\frac{1}{2}H^2(P_{\theta^{(j)}},P_{\theta})$   combining  \eqref{pr:lower:e3} with 
 \eqref{pr:lower:e4} yields
\begin{equation*}
\Bigl\{\Ex_{{\theta^{(j)}}}|\skalarV{\widetilde{\beta}-\beta_{\theta^{(j)}},\psi_j}|^2+ \Ex_{{\theta}}|\skalarV{\widetilde{\beta}-\beta_{\theta},\psi_j}|^2\Bigr\}\geqslant\frac{1}{2} b_j^2,\quad j=1,\dotsc, \kstar.
 \end{equation*}
From this  we conclude for each estimator $\widetilde\beta$ that
\begin{align*}
\sup_{\beta \in \cF_\bw^\br} &\Ex\normV{\widetilde\beta -\beta}_\hw^2 \geqslant \sup_{\theta\in \{-1,1\}^{\kstar}} \Ex_\theta\normV{\widetilde\beta -\beta_{\theta}}_\hw^2\\
&\geqslant \frac{1}{2^{\kstar}}\sum_{\theta\in \{-1,1\}^{\kstar}}\sum_{j=1}^{\kstar}\hw_j\,\Ex_{{\theta}}|\skalarV{\widetilde{\beta}-\beta_{\theta},\psi_j}|^2\\
&= \frac{1}{2^{\kstar}}\sum_{\theta\in \{-1,1\}^{\kstar}}\frac{1}{2}\sum_{j=1}^{\kstar}\hw_j\Bigl\{\Ex_{{\theta}}|\skalarV{\widetilde{\beta}-\beta_{\theta},\psi_j}|^2+\Ex_{{\theta^{(j)}}}|\skalarV{\widetilde{\beta}-\beta_{\theta^{(j)}},\psi_j}|^2 \Bigr\}\\
&\geqslant\frac{1}{4} \sum_{j=1}^{\kstar}  b_j^2\, \hw_j \geqslant \frac{1}{4}\min \left( \frac{\sigma^2}{2\,d},  \frac{\rho}{\triangle}\right) \,\frac{\dstar }{\triangle},
\end{align*}
where the last inequality follows again from  \eqref{app:l1:lower:e1} in Lemma \ref{app:l1:lower},  which completes the proof.\hfill$\square$\\
\paragraph{Proof of the upper bound.}
\begin{proof}[\noindent\textcolor{darkred}{\sc Proof of Theorem \ref{res:upper}.}]The proof is based on the  decomposition \eqref{res:dec}, where  we show  below under the condition $X\in \cX_{\eta}^{8}$, $\epsilon\in \cE_\eta^{8}$ and $\lambda_j\geqslant 2\alpha$, $1\leqslant j\leqslant m$, for some generic constant $C>0$ the following two bounds
\begin{gather}\label{pr:res:upper:gen:e1}
\Ex\normV{\widehat{\beta}- \widetilde{\beta}_m}^2_\hw
\leqslant C \sum_{j=1}^m \frac{\hw_j}{n\, \lambda_j} \, \Bigl\{ \normV{\beta}^{2}\, \Ex\normV{X}^2+\sigma^{2}\Bigr\}\, \eta \,\Bigl\{  \lambda_j^2/(n\alpha)^2  + (1/n) +1\Bigr\},
\\\label{pr:res:upper:gen:e2}
\Ex\normV{\widetilde{\beta}_m-\beta}^2_\hw\leqslant C \{\hw_m/\bw_m + \eta/n \}\,\normV{\beta}_\bw^2,
\end{gather}
Consequently, for all $\beta\in \cF_\bw^\br$ and $(\lambda_j)_{j\geqslant 1} \in \cS_\lw^\ld$, i.e., $\lambda_j\leqslant \ld\,\lw_j\leqslant d $ and  $\Ex\normV{X}^2\leqslant \ld \Lw$, follows
\begin{equation*}
\Ex\normV{\widehat{\beta}- \beta}^2_\hw\leqslant C \,  \Bigl\{  d\,(\ld^2/(n\alpha)^2 + 1/n +1)\, \sum_{j=1}^m \frac{\hw_j}{n\, \lw_j}  + \hw_m/\bw_m + 1/n \Bigr\}\,\eta\,[ \br\, d\,  \Lw+\sigma^2].
\end{equation*}
Let $\kstar$ and $\dstar$ be given by \eqref{res:lower:def:md} for some $\triangle\geqslant 1$ then the condition on $m$ and   $\alpha$,  i.e.,  $m=\kstar$ and  $\alpha = (1/n) \min(1,\gamma_{\kstar}/(2\ld\triangle))$, implies 
\begin{equation*}
\Ex\normV{\widehat{\beta}- \beta}^2_\hw\leqslant C \,   d\,(\ld^2/(n\alpha)^2 + 1/n +1)\,\triangle\,\max(\dstar,1/n)\,\eta\,[ \br\, d\,  \Lw+\sigma^2].
\end{equation*}
because $\hw_m/\bw_m=\dstar$, $\sum_{j=1}^m {\hw_j}/(n\, \lw_j)\leqslant \triangle \dstar$ and  $\lambda_j\geqslant 2\alpha$, $1\leqslant j\leqslant m$ by using that $\lw_{\kstar} \geqslant \bw_{\kstar}/(n\triangle)$ and $(\lambda_j)_{j\geqslant 1} \in \cS_\lw^\ld$. Hence,  from $n\alpha \geqslant 1/(2\ld\triangle)$ follows the result.

Proof of (\ref{pr:res:upper:gen:e1}). By using  $T_{n,j}$ introduced in \eqref{app:l:def} we  obtain the identity
 \begin{equation}\label{pr:res:upper:gen:e1:1}
\Ex\normV{\widehat{\beta}_\hw- \widetilde{\beta}^\alpha_{\hw}}^2_\hw
=  \sum_{j=1}^m \frac{\hw_j}{\lambda_j} \cdot \Ex\Bigl[\frac{T_{n,j}^2}{\lambda_j}\cdot |\lambda_j/
\widehat{\lambda}_j|^2 \1\{\widehat{\lambda}_j\geqslant \alpha\}\Bigr].
 \end{equation}
By using the elementary inequality $1/2\leqslant |\widehat{\lambda}_j/\lambda_j-1|^2 + |\widehat{\lambda}_j/\lambda_j|^2$ it follows that 
\begin{equation*}|\lambda_j/
\widehat{\lambda}_j|^2 \1\{\widehat{\lambda}_j\geqslant \alpha\}\leqslant 2\Bigl\{ 2\, (\lambda_j/\alpha)^2\,  |\widehat\lambda_j/\lambda_j-1|^4 + 2 |\widehat\lambda_j/\lambda_j-1|^2+1\Bigr\}.\end{equation*}
Therefore, by combination of  the last estimate and \eqref{pr:res:upper:gen:e1:1} we have
 \begin{multline*}
\Ex\normV{\widehat{\beta}- \widetilde{\beta}^\alpha_{\hw}}^2_\hw
\leqslant   4 \sum_{j=1}^m \frac{\hw_j}{\lambda_j} \cdot \Bigl(\Ex|T_{n,j}|^4\Bigr)^{1/2}\cdot \Bigl\{  (\lambda_j/\alpha)^2 \Bigl(\Ex|\widehat{\lambda}_j/\lambda_j-1|^8\Bigr)^{1/2} + \Bigl(\Ex|\widehat{\lambda}_j/\lambda_j-1|^4\Bigr)^{1/2} +1\Bigr\}
 \end{multline*}
The estimate (\ref{pr:res:upper:gen:e1}) follows now from \eqref{app:l1:upp:e1} and \eqref{app:l1:upp:e3} in Lemma \ref{app:l1:upp}.

Proof of (\ref{pr:res:upper:gen:e2}). Following  along the lines of the proof of (\ref{pr:res:prop:cons:e2})  we obtain
\begin{equation*}\Ex\normV{\widetilde{\beta}_m-\beta}_\hw^2\leqslant 2\{ \normV{\beta_m-\beta}^2_\hw+ C (\eta/n)\, \normV{\beta_m}_\hw^2\},
\end{equation*}
where under the condition $\lambda_j\geqslant 2\alpha$ for each $1\leqslant j\leqslant m$ we have used that $P(\widehat{\lambda}_j<\alpha) \leqslant C\,\eta/n$. Then, under Assumption \ref{ass:reg}, i.e., $(\hw_j/\bw_j)$ is non-increasing, the usual estimate $\normV{\beta_m-\beta}^2_\hw\leqslant  \hw_m/\bw_m\normV{\beta}_\bw^2$ implies (\ref{pr:res:upper:gen:e2}), which completes the proof. \end{proof}
\paragraph{Technical assertions.}\hfill\\[1ex]
The following two lemma gather technical results used in the proof of Proposition \ref{res:prop:cons}, Theorem \ref{res:lower} and Theorem \ref{res:upper}.
\begin{lem}\label{app:l1:upp} Suppose $X\in\cX_{\eta}^{4m}$ and $\epsilon\in\cE^{4m}_\eta$, $m\in\N$. Then for some constant $C>0$ only depending on $m$ we have
\begin{gather}\label{app:l1:upp:e1}
\sup_{j\in\N}\Bigl\{\lambda_j^{-m}\cdot  \Ex |T_{n,j}|^{2m}\Bigr\} \leq C \cdot n^{-m}\cdot\{ \normV{\beta}^{2m}\cdot (\Ex\normV{X}^2)^m+\sigma^{2m}\}\cdot \eta,\\
\label{app:l1:upp:e2}
\sup_{j\in\N}\Ex |\widehat{\lambda}_j/\lambda_j-1|^{2m}\leq C\cdot n^{-m}\cdot \eta.
\end{gather}
If in addition $w_1\geq 2$ and $w_2\leq 1/2$, then we obtain
\begin{gather}\label{app:l1:upp:e3}
\sup_{j\in\N}P(\widehat{\lambda}_j/\lambda_j\geq w_1)\leq C \cdot n^{-m}\cdot \eta\;\text{ and }\;
\sup_{j\in\N}P(\widehat{\lambda}_j/\lambda_j< w_2)\leq C \cdot n^{-m}\cdot \eta.
\end{gather}
\end{lem}
\noindent\textcolor{darkred}{\sc Proof.} Let $\zeta_{ij}:=\sum_{l\ne j }\beta_l {X}_{il}$, $i=1,\dotsc,n$ and $j\in\N$. Then we have
\begin{equation*}T_{n,j}=\frac{1}{n}\sum_{i=1}^n\{ \zeta_{ij}+ \sigma \epsilon_i\} X_{ij}=:T_1+T_2,\end{equation*}
where we bound below each summand separately, that is
\begin{align}\label{app:l1:upp:e1:1}
\Ex|T_1|^{2m}&\leq C\cdot \frac{\lambda_j^m}{n^{m}}\cdot  \normV{\beta}^{2m} \cdot  (\Ex\normV{X}^2)^m\cdot \eta,\\\label{app:l1:upp:e1:2}
\Ex|T_2|^{2m}&\leq C\cdot \frac{\lambda_j^m}{n^{m}}\cdot\sigma^{2m} \cdot \eta
\end{align}
for some $C>0$ only depending on $m$. Consequently, the  inequality \eqref{app:l1:upp:e1} follows from \eqref{app:l1:upp:e1:1} and \eqref{app:l1:upp:e1:2}. Consider $T_1$. For each $j\in \N$ the random variables $(\zeta_{ij}\cdot X_{ij})$,  $i=1,\dots,n,$ are independent and identically distributed with mean zero. From Theorem 2.10 in \cite{Petrov1995} we conclude $\Ex |T_{1}|^{2m} \leq C n^{-m}\Ex |\zeta_{1j} X_{1j}|^{2m}$ for some  constant $C>0$ only depending on $m$. Then  we claim that \eqref{app:l1:upp:e1:1} follows in case of  $T_1$ from  the Cauchy-Schwarz inequality together with   $X_1\in\cX_{\eta}^{4m}$, i.e.,  $\sup_{j} \Ex|X_{1j}/\sqrt{\lambda_j}|^{4m}\leq \eta$. Indeed, we have
$$\Ex |\zeta_{1j}X_{1j}|^{2m}
\leq (\sum_{l\ne j}\beta_l^2)^m\sum_{l_1\neq j}\dots\sum_{l_m\ne j} \Ex |X_{1j}|^{2m}\prod_{k=1}^m |X_{1l_k}|^2\leq \normV{\beta}^{2m}\cdot \lambda_j^m \cdot  (\sum_{l\neq j}\lambda_{l})^m\cdot \eta.$$
Consider $T_2$. \eqref{app:l1:upp:e1:2} follows in analogy to the case of $T_1$, because  $\{\sigma\, \epsilon_{i}\, X_{ij} \}$ are independent and identically distributed with mean zero, and 
$\Ex |\sigma \cdot \epsilon_1\cdot X_{1j}|^{2m}\leq \sigma^{2m} \cdot \lambda_j^m \cdot \eta$.

Proof of (\ref{app:l1:upp:e2}).  Since $\{(|X_{ij}|^2/\lambda_j-1)\}$  are independent and identically distributed with mean zero, and 
$\Ex| X_{1j}^2/\lambda_j|^{2m}\leq \eta$, the result follows by applying Theorem 2.10 in \cite{Petrov1995}.

Proof of (\ref{app:l1:upp:e3}). If $w\geq 2$ then  $P(\widehat{\lambda}_j/\lambda_j\geq w)\leq P(|\widehat{\lambda}_j/\lambda_j-1|\geq 1)$. Thus  applying Markov's inequality together with \eqref{app:l1:upp:e2} implies the first bound in (\ref{app:l1:upp:e3}), while the second follows in analogy,  which proves the lemma.\hfill$\square$\\

\begin{lem}\label{app:l1:lower}Let $\kstar\in\N$ and $\dstar$ be chosen such that  \eqref{res:lower:def:md} is satisfied for some $\triangle\geqslant 1$.  Consider a (infinite) vector $b$ with components $b_j$ satisfying  \begin{equation}\label{app:l1:lower:b}b_j^2= \frac{\zeta}{n\cdot \lw_j},\quad j\in\N, \quad \text{ with }\quad \zeta:=\min \left( \sigma^2/(2d),  \br/\triangle\right),\end{equation}
then we have for all $j\in\N$
\begin{gather}\label{app:l1:lower:e1}
\frac{2\,d\, n }{\sigma^2}b^2_j \,\lw_j\leqslant 1,\quad
\sum_{j=1}^{\kstar}b^2_j \, \bw_j\leqslant \rho,\quad \mbox{and}\quad
 \sum_{j=1}^{\kstar}\,b^2_j\,  \hw_j\geqslant \min \left( \frac{\sigma^2}{2{d}},  \frac{\rho}{\triangle}\right) \,\frac{\max(\dstar,1/n)}{\triangle}.
\end{gather}
\end{lem}
\noindent\textcolor{darkred}{\sc Proof.} The first inequality in \eqref{app:l1:lower:e1}  follows trivially by using the definition of $\zeta$. Since by Assumption \ref{ass:reg} the sequence $(\bw_j/\hw_j)$ is nondecreasing  the definition of $\kstar$ given in \eqref{res:lower:def:md} implies the second estimate in \eqref{app:l1:lower:e1}, i.e.,  $\sum_{j=1}^{\kstar}b^2_{j} \bw_j \leqslant 
\zeta (\bw_{\kstar}/\hw_{\kstar})\sum_{j=1}^{\kstar}\hw_j/(n \lw_j) \leqslant\zeta\triangle\leqslant \rho$. To deduce the third inequality in \eqref{app:l1:lower:e1}  from the definition of $\kstar$ and $\dstar$ observe that $ \sum_{j=1}^{\kstar}b^2_j \hw_j =  \dstar \,\zeta \,(\bw_{\kstar}/\hw_{\kstar}) \sum_{j=1}^{\kstar}\hw_{j}/{(n\, \lw_j)} \geqslant \dstar \,\zeta/\triangle$ and $ \sum_{j=1}^{\kstar}b^2_j \hw_j \geqslant  \zeta /n$ since $\hw_1/\lw_1=1$, which proves the lemma.\hfill$\square$
\subsection{Proofs of Section \ref{sec:ex}}\label{app:proofs:ex}
\paragraph{The mean prediction error.}
\begin{proof}[\noindent\textcolor{darkred}{\sc Proof of Proposition \ref{MPE:lower}.}] Given the eigenvalues $(\lambda_j)$ of $\op$ satisfy a link condition, that is $(\lambda_j) \in \cS_{\lw}^{\ld}$, $\ld\geqslant 1$.   It follows  that $\Ex \normV{\widehat\beta-\beta}_\lw^2 \asymp_d \Ex\skalarV{\op(\widehat\beta-\beta),\widehat\beta-\beta}$. Therefore, we can apply the general results  by considering the $\cF_\hw$-risk with $\hw\equiv\lw$.  Furthermore, in case (i) the definition of $\bw \equiv w^p$  and $\lw$ imply together $(\bw_{\kstar}/\hw_{\kstar})\sum_{j=1}^{\kstar}  \hw_j/\lw_j =  \kstar^{2a+2p+1}$. It follows that the condition on $\kstar$ and $\dstar$ given in \eqref{res:lower:def:md} of Theorem \ref{res:lower} can be rewritten as $\kstar\sim n^{1/(2p+2a+1)}$ and $\dstar\sim n^{-(2p+2a)/(2p+2a+1)}$, respectively. On the other hand, in case (ii)   $(\bw_{\kstar}/\hw_{\kstar})\sum_{j=1}^{\kstar}  \hw_j/\lw_j= \kstar^{2p+1}\exp(\kstar^{2a})$ implies that the condition on $\kstar$ and $\dstar$ writes  $\kstar\sim (\log n)^{1/(2a)}$ and $\dstar\sim n^{-1}(\log n)^{1/(2a)}$, respectively. Consequently, the lower bounds in Proposition \ref{MPE:lower} follow by applying Theorem \ref{res:lower}.\end{proof}

\begin{proof}[\noindent\textcolor{darkred}{\sc Proof of Proposition \ref{MPE:upper}.}] Since in both cases  the condition on the dimension parameter $m$ and the threshold $\alpha$ ensures that  $m\sim \kstar$  and $\alpha \sim 1/n$ (see the proof of Proposition \ref{MPE:lower}) the result follows from Theorem \ref{res:upper}.\end{proof}

\paragraph{The estimation of derivatives.}
\begin{proof}[\noindent\textcolor{darkred}{\sc Proof of Proposition \ref{MSE:lower}.}] Due to $\Ex\normV{\widetilde\beta^{(s)}-\beta^{(s)}}^2 \asymp_{(2\pi)^{2s}} \Ex \normV{\widetilde\beta-\beta}_{w^s}^2$, $0\leqslant s \leqslant p$,  we can apply again the general results  by considering the $\cF_\hw$-risk with $\hw\equiv w^s$. In case (i) the well-known  approximation $\sum_{j=1}^{m} j^{r}\sim m^{r+1}$ for $r>0$ together with the definition of $\bw \equiv w^p$  and $\lw$  implies  $(\bw_{\kstar}/\hw_{\kstar})\sum_{j=1}^{\kstar}  \hw_j/\lw_j \sim  \kstar^{2a+2p+1}$. It follows that the condition on $\kstar$ and $\dstar$ given in \eqref{res:lower:def:md} of Theorem \ref{res:lower} writes $\kstar\sim n^{1/(2p+2a+1)}$ and $\dstar\sim n^{-(2p-2s)/(2p+2a+1)}$, respectively. On the other hand, in case (ii) by applying Laplace's Method (c.f. chapter 3.7 in \cite{Olver1974}) the definition of $\bw \equiv w^p$  and $\lw$  
imply  $(\bw_{\kstar}/\hw_{\kstar})\sum_{j=1}^{\kstar}  \hw_j/\lw_j\sim  \kstar^{2p}\exp(\kstar^{2a})$. Therefore,  the condition on $\kstar$ and $\dstar$ can be rewritten as   $\kstar\sim (\log n)^{1/(2a)}$ and $\dstar\sim n^{-1}(\log n)^{1/(2a)}$, respectively. Consequently, the lower bounds in Proposition \ref{MPE:lower} follow by applying Theorem \ref{res:lower}.\end{proof}

\begin{proof}[\noindent\textcolor{darkred}{\sc Proof of Proposition \ref{MSE:upper}.}]Since in both cases  the condition on the dimension parameter $m$ and the threshold $\alpha$ ensures that  $m\sim \kstar$  and $\alpha \sim 1/n$ (see the proof of Proposition \ref{MSE:lower}) the result follows from Theorem \ref{res:upper}.\end{proof}

\bibliography{BiB-NP-FLM-SOR}
\end{document}